\documentclass[11pt,reqno]{amsart}

\usepackage{amsmath, amsthm, amsopn, amssymb, enumerate, fullpage}

\usepackage{latexsym}
\usepackage{graphics}
\usepackage{graphicx}
\usepackage{verbatim} 
\usepackage{hyperref} 
\usepackage[mathscr]{eucal}

\newtheorem{thm}{Theorem}[section]
\newtheorem{conj}[thm]{Conjecture}
\newtheorem{lemma}[thm]{Lemma}
\newtheorem{cor}[thm]{Corollary}
\newtheorem*{HCL}{The hypergraph container lemma}
\newtheorem*{HCL3}{The hypergraph container lemma for 3-uniform hypergraphs}

\newcommand{\Bin}{\textup{Bin}}
\newcommand{\col}{\mathrm{col}}
\newcommand{\eps}{\varepsilon}
\newcommand{\ex}{\mathrm{ex}}
\newcommand{\Fnind}{\mathcal{F}_n^{\textup{ind}}}
\newcommand{\Gnpind}{G _{n,p}^{\mathrm{ind}}}
\newcommand{\vol}{\mathrm{vol}}

\newcommand{\scolon}{\,\colon}

\def\A{\mathcal{A}}
\def\C{\mathcal{C}}
\def\cF{\mathcal{F}}
\def\M{\mathcal{M}}
\def\G{\mathcal{G}}
\def\I{\mathcal{I}}
\def\cH{\mathcal{H}}
\def\N{\mathbb{N}}
\def\Z{\mathbb{Z}}

\def\R{\mathbb{R}}
\def\Z{\mathbb{Z}}
\def\P{\mathcal{P}}
\def\Pr{\mathbb{P}}
\def\Ex{\mathbb{E}}
\def\S{\mathcal{S}}

\def\le{\leqslant}
\def\ge{\geqslant}

\makeatletter
\def\Ddots{\mathinner{\mkern1mu\raise\p@
\vbox{\kern7\p@\Hbbox{.}}\mkern2mu
\raise4\p@\Hbbox{.}\mkern2mu\raise7\p@\Hbbox{.}\mkern1mu}}
\makeatother

\begin{document}

\title{The method of hypergraph containers}

\author{J\'ozsef Balogh}
\address{Department of Mathematics, University of Illinois at Urbana-Champaign, Urbana, Illinois 61801, USA}
\email{jobal@math.uiuc.edu}

\author{Robert Morris}
\address{IMPA, Estrada Dona Castorina 110, Jardim Bot\^anico, Rio de Janeiro, 22460-320, Brazil}
\email{rob@impa.br}

\author{Wojciech Samotij}
\address{School of Mathematical Sciences, Tel Aviv University, Tel Aviv 6997801, Israel}
\email{samotij@post.tau.ac.il}

\thanks{JB is partially supported by NSF Grant DMS-1500121 and by the Langan Scholar Fund (UIUC); RM is partially supported by CNPq (Proc.~303275/2013-8), by FAPERJ (Proc.~201.598/2014), and by ERC Starting Grant 680275 MALIG; WS is partially supported by the Israel Science Foundation grant 1147/14.}

\begin{abstract}
  In this survey we describe a recently-developed technique for bounding the number (and controlling the typical structure) of finite objects with forbidden substructures. This technique exploits a subtle clustering phenomenon exhibited by the independent sets of uniform hypergraphs whose edges are sufficiently evenly distributed; more precisely, it provides a relatively small family of `containers' for the independent sets, each of which contains few edges. We attempt to convey to the reader a general high-level overview of the method, focusing on a small number of illustrative applications in areas such as extremal graph theory, Ramsey theory, additive combinatorics, and discrete geometry, and avoiding technical details as much as possible. 
\end{abstract}

\maketitle

\section{Introduction}

Numerous well-studied problems in combinatorics concern families of discrete objects which avoid certain forbidden configurations, such as the family of $H$-free graphs\footnote{A graph is $H$-free if it does not contain a subgraph isomorphic to $H$.} or the family of sets of integers containing no $k$-term arithmetic progression. The most classical questions about these families relate to the size and structure of the extremal examples; for example, Tur\'an~\cite{T41} determined the unique $K_r$-free graph on $n$ vertices with the most edges and Szemer\'edi~\cite{Sz75} proved that every set of integers of positive upper density contains arbitrarily long arithmetic progressions. In recent decades, partly motivated by applications to areas such as Ramsey theory and statistical physics, there has been increasing interest in problems relating to the typical structure of a (e.g., uniformly chosen) member of one of these families and to extremal questions in (sparse) random graphs and random sets of integers. Significant early developments in this direction include the seminal results obtained by Erd\H{o}s, Kleitman, and Rothschild~\cite{EKR}, who proved that almost all triangle-free graphs are bipartite, by Kleitman and Winston~\cite{KW82}, who proved that there are $2^{\Theta(n^{3/2})}$ $C_4$-free graphs on $n$ vertices, and by Frankl and R\"odl~\cite{FR}, who proved that if $p \gg 1/\sqrt{n}$, then with high probability every 2-colouring of the edges of $G(n,p)$ contains a monochromatic triangle.

An important recent development in this area was the discovery that, perhaps surprisingly, it is beneficial to consider such problems in the more abstract (and significantly more general) setting of independent sets in hypergraphs. This approach was taken with stunning success by Conlon and Gowers~\cite{CG}, Friedgut, R\"odl, and Schacht~\cite{FRS}, and Schacht~\cite{Sch} in their breakthrough papers on extremal and Ramsey-type results in sparse random sets. To give just one example of the many important conjectures resolved by their work, let us consider the random variable
\[
  \ex\big( G(n,p), H \big) = \max\big\{ e(G) \scolon  H \not\subset G \subset G(n,p) \big\},
\]
which was first studied (in the case $H = K_3$) by Frankl and R\"odl~\cite{FR}. The following theorem was conjectured by Haxell, Kohayakawa, and \L uczak~\cite{HKL-odd, HKL-even} and proved (independently) by Conlon and Gowers~\cite{CG} and by Schacht~\cite{Sch}.

\begin{thm}
  \label{thm:Turan:Gnp}
  Let $H$ be a graph with at least two edges and suppose that $p \gg n^{-1/m_2(H)}$, where $m_2(H)$ is the so-called $2$-density\footnote{To be precise, $m_2(H) = \max\big\{ \frac{e(F) - 1}{v(F) - 2} \scolon F \subset H, \, v(F) \ge 3 \big\}$.} of $H$. Then
  \begin{equation*}\label{eq:thm:Turan:Gnp}
    \ex\big( G(n,p), H \big) = \bigg( 1 - \frac{1}{\chi(H) - 1} + o(1) \bigg) p {n \choose 2}
  \end{equation*}
  asymptotically almost surely (a.a.s.), that is, with probability tending to $1$ as $n \to \infty$. 
\end{thm}

It is not hard to show that $\ex\big( G(n,p), H \big) = \big( 1 + o(1) \big) p {n \choose 2}$ a.a.s.\ if $n^{-2} \ll p \ll n^{-1/m_2(H)}$ and so the assumption on $p$ in Theorem~\ref{thm:Turan:Gnp} is optimal. We remark that in the case when $H$ is a clique even more precise results are known, due to work of DeMarco and Kahn~\cite{DK2, DK1}, who proved that if $p \gg n^{-1/m_2(H)} (\log n)^{2/(r+1)(r - 2)}$, then with high probability the largest $K_{r+1}$-free subgraph of $G(n,p)$ is $r$-partite, which is again essentially best possible. We refer the reader to an excellent recent survey of R\"odl and Schacht~\cite{RS} for more details on extremal results in sparse random sets.

In this survey we will describe an alternative approach to the problem of understanding the family of independent sets in a hypergraph, whose development was inspired by the work in~\cite{CG,FRS,Sch} and also strongly influenced by that of Kleitman and Winston~\cite{KW82} and Sapozhenko~\cite{Sap01,Sap03,Sap05}. This technique, which was developed independently by the authors of this survey~\cite{BMS} and by Saxton and Thomason~\cite{ST}, has turned out to be surprisingly powerful and flexible. It allows one to prove enumerative, structural, and extremal results (such as Theorem~\ref{thm:Turan:Gnp}) in a wide variety of settings. It is known as the \emph{method of hypergraph containers}.

To understand the essence of the container method, it is perhaps useful to consider as an illustrative example the family $\cF_n(K_3)$ of triangle-free graphs on (a given set of) $n$ vertices. Note that the number of such graphs is at least $2^{\lfloor n^2/4 \rfloor}$, since every bipartite graph is triangle-free.\footnote{In particular, every subgraph of the complete bipartite graph with $n$ vertices and $\lfloor n^2/4 \rfloor$ edges is triangle-free.} However, it turns out that there exists a vastly smaller family $\G_n$ of graphs on $n$ vertices, of size $n^{O(n^{3/2})}$, that forms a set of \emph{containers} for $\cF_n(K_3)$, which means that for every $H \in \cF_n(K_3)$, there exists a $G \in \G_n$ such that $H \subset G$. A remarkable property of this family of containers is that each graph $G \in \G_n$ is `almost triangle-free' in the sense that it contains `few' triangles. It is not difficult to use this family of containers, together with a suitable `supersaturation' theorem, to prove Theorem~\ref{thm:Turan:Gnp} in the case $H = K_3$ or to show, using a suitable `stability' theorem, that almost all triangle-free graphs are `almost bipartite'. We will discuss these two properties of the family of triangle-free graphs in much more detail in Section~\ref{basic:sec}.

In order to generalize this container theorem for triangle-free graphs, it is useful to first restate it in the language of hypergraphs. To do so, consider the 3-uniform hypergraph $\cH$ with vertex set $V(\cH) = E(K_n)$ and edge set
\[
  E(\cH) = \big\{ \{e_1,e_2,e_3\} \subset E(K_n) \scolon \textup{$e_1, e_2, e_3$ form a triangle} \big\}.
\]
We shall refer to $\cH$ as the `hypergraph that encodes triangles' and emphasize that (somewhat confusingly) the vertices of this hypergraph are the edges of the complete graph $K_n$. Note that $\cF_n(K_3)$ is precisely the family $\I(\cH)$ of independent sets of $\cH$, so we may rephrase our container theorem for triangle-free graphs as follows:
\begin{align*}
& \text{``There exists a relatively small family $\C$ of subsets of $V(\cH)$, each containing only few}\\
& \text{edges of $\cH$, such that every independent set $I \in \I(\cH)$ is contained in some member of $\C$."}
\end{align*}
There is nothing special about the fine structure of the hypergraph encoding triangles that makes the above statement true. On the contrary, the method of containers allows one to prove that a similar phenomenon holds for a large class of $k$-uniform hypergraphs, for each $k \in \N$. In the case $k=3$, a sufficient condition is the following assumption on the distribution of the edges of a~$3$-uniform hypergraph $\cH$ with average degree $d$: each vertex of $\cH$ has degree at most $O(d)$ and each pair of vertices lies in at most $O(\sqrt{d})$ edges of $\cH$. For the hypergraph that encodes triangles, both conditions are easily satisfied, since each edge of $K_n$ is contained in exactly $n-2$ triangles and each pair of edges is contained in at most one triangle. The conclusion of the container lemma (see Sections~\ref{basic:sec} and~\ref{keylemma:sec}) is that each independent set $I$ in a $3$-uniform hypergraph $\cH$ satisfying these conditions has a \emph{fingerprint} $S \subset I$ of size $O\big( v(\cH) / \sqrt{d} \big)$ that is associated with a set $X(S)$ of size $\Omega\big( v(\cH) \big)$ which is disjoint from $I$. The crucial point is that the set $X(S)$ depends only on $S$ (and not on $I$) and therefore the number of sets $X(S)$ is bounded from above by the number of subsets of the vertex set $V(\cH)$ of size $O\big( v(\cH) / \sqrt{d} \big)$. In particular, each independent set of $\cH$ is contained in one of at most $v(\cH)^{O(v(\cH)/\sqrt{d})}$ sets of size at most $(1 - \delta)v(\cH)$, for some constant $\delta > 0$. By iterating this process, that is, by applying the container lemma repeatedly to the subhypergraphs induced by the containers obtained in earlier applications, one can easily prove the container theorem for triangle-free graphs stated (informally) above.

Although the hypergraph container lemma (see Section~\ref{keylemma:sec}) was discovered only recently (see~\cite{BMS,ST}), several theorems of the same flavour (though often in very specific settings) appeared much earlier in the literature. The earliest container-type argument of which we are aware appeared (implicitly) over 35 years ago in the work of Kleitman and Winston on bounding the number of lattices~\cite{KW80} and of $C_4$-free graphs~\cite{KW82}, which already contained some of the key ideas needed for the proof in the general setting; see~\cite{Sa15} for details. Nevertheless, it was not until almost 20 years later that Sapozhenko~\cite{Sap01,Sap03,Sap05} made a systematic study of containers for independent sets in graphs (and coined the name \emph{containers}). Around the same time, Green and Ruzsa~\cite{GrRu04} obtained (using Fourier analysis) a container theorem for sum-free subsets of $\Z/p\Z$.

More recently, Balogh and Samotij~\cite{BSmm,BSst} generalized the method of~\cite{KW82} to count $K_{s,t}$-free graphs, using what could be considered to be the first container theorem for hypergraphs of uniformity larger than two. Finally, Alon, Balogh, Morris, and Samotij~\cite{ABMS1,ABMS2} proved a general container theorem for 3-uniform hypergraphs and used it to prove a sparse analogue of the Cameron--Erd\H{o}s conjecture. Around the same time, Saxton and Thomason~\cite{SaTh12} developed a simpler version of the method and applied it to the problem of bounding the list chromatic number of hypergraphs. In particular, the articles~\cite{ABMS1} and~\cite{SaTh12} can be seen as direct predecessors of~\cite{BMS} and~\cite{ST}.

The rest of this survey is organised as follows. In Section~\ref{basic:sec}, we warm up by stating a container lemma for 3-uniform hypergraphs, giving three simple applications to problems involving triangle-free graphs and a more advanced application to a problem in discrete geometry that was discovered recently by Balogh and Solymosi~\cite{BS}. Next, in Section~\ref{keylemma:sec}, we state the main container lemma and provide some additional motivation and discussion of the statement and in Section~\ref{counting:Hfree:sec} we describe an application to counting $H$-free graphs. Finally, in Sections~\ref{sec:many-colours}--\ref{more:applications:sec}, we state and discuss a number of additional applications, including to multi-coloured structures (e.g., metric spaces), asymmetric structures (e.g., sparse members of a hereditary property), hypergraphs of unbounded uniformity (e.g., induced Ramsey numbers, $\eps$-nets), number-theoretic structures (e.g., Sidon sets, sum-free sets, sets containing no $k$-term arithmetic progression), sharp thresholds in Ramsey theory, and probabilistic embedding in sparse graphs.

\section{Basic applications of the method}\label{basic:sec}

In this section we will provide the reader with a gentle introduction to the container method, focusing again on the family of triangle-free graphs. In particular, we will state a version of the container lemma for 3-uniform hypergraphs and explain (without giving full details) how to deduce from it bounds on the largest size of a triangle-free subgraph of the random graph $G(n,p)$, statements about the typical structure of a (sparse) triangle-free graph, and how to prove that every $r$-colouring of the edges of $G(n,p)$ contains a monochromatic triangle. To give a simple demonstration of the flexibility of the method, we will also describe a slightly more complicated application to a problem in discrete geometry. 


In order to state the container lemma, we need a little notation. Given a hypergraph $\cH$, let us write $\Delta_\ell(\cH)$ for the maximum degree of a set of $\ell$ vertices of $\cH$, that is,
\[
  \Delta_\ell(\cH) = \max\big\{ d_\cH(A) \scolon A \subset V(\cH), \, |A| = \ell \big\},
\]
where $d_\cH(A) = \big| \big\{ B \in E(\cH) \scolon A \subset B \big\} \big|$, and $\I(\cH)$ for the collection of independent sets of $\cH$.

\begin{HCL3}
  For every $c > 0$, there exists $\delta > 0$ such that the following holds. Let $\mathcal{H}$ be a $3$-uniform hypergraph with average degree $d \ge \delta^{-1}$ and suppose that
  \[
    \Delta_1(\cH) \le c \cdot d \qquad \text{and} \qquad \Delta_2(\cH) \le c \cdot \sqrt{d}.
  \]
  Then there exists a collection $\C$ of subsets of $V(\mathcal{H})$ with
  \[
    |\C| \le \binom{v(\cH)}{v(\cH)/\sqrt{d}}
  \]
  such that
  \begin{enumerate}
  \item[$(a)$] for every $I \in \I(\cH)$, there exists $C \in \C$ such that $I \subset C$,
  \item[$(b)$] $|C| \le (1 - \delta) v(\cH)$ for every $C \in \C$.
  \end{enumerate}
\end{HCL3}

In order to help us understand the statement of this lemma, let us apply it to the hypergraph $\cH$ that encodes triangles in $K_n$, defined in the Introduction. Recall that this hypergraph satisfies
\[
  v(\cH) = {n \choose 2}, \qquad \Delta_2(\cH) = 1, \qquad \text{and} \qquad d_\cH(v) = n - 2
\]
for every $v \in V(\cH)$. We may therefore apply the container lemma to $\cH$, with $c = 1$, to obtain a collection $\C$ of $n^{O(n^{3/2})}$ subsets of $E(K_n)$ (that is, graphs on $n$ vertices) with the following properties:
\begin{enumerate}
\item[$(a)$] Every triangle-free graph is a subgraph of some $C \in \C$.
\item[$(b)$] Each $C \in \C$ has at most $(1 - \delta) e(K_n)$ edges.
\end{enumerate}

Now, if there exists a container $C \in \C$ with at least $\eps n^3$ triangles, then take each such $C$ and apply the container lemma to the subhypergraph $\cH[C]$ of $\cH$ induced by $C$, i.e., the hypergraph that encodes triangles in the graph $C$. Note that the average degree of $\cH[C]$ is at least $6\eps n$, since each triangle in $C$ corresponds to an edge of $\cH[C]$ and $v(\cH[C]) = |C| \le e(K_n)$. Since (trivially) $\Delta_\ell(\cH[C]) \le \Delta_\ell(\cH)$, it follows that we can apply the lemma with $c = 1/\eps$ and replace $C$ by the collection of containers for $\I(\cH[C])$ given by the lemma. 

Let us iterate this process until we obtain a collection $\C$ of containers, each of which has fewer than $\eps n^3$ triangles. How large is the final family $\C$ that we obtain? Note that we apply the lemma only to hypergraphs with at most $\binom{n}{2}$ vertices and average degree at least $6\eps n$ and therefore produce at most $n^{O(n^{3/2})}$ new containers in each application, where the implicit constant depends only on~$\eps$. Moreover, each application of the lemma shrinks a container by a factor of $1 - \delta$, so after a bounded (depending on $\eps$) number of iterations every container will have fewer than $\eps n^3$ triangles (since $\Delta_1(\cH) < n$, then every graph with at most $\eps n^2$ edges contains fewer than $\eps n^3$ triangles).

The above argument yields the following container theorem for triangle-free graphs.

\begin{thm}\label{thm:CT:triangles}
  For each $\eps > 0$, there exists $C > 0$ such that the following holds. For each $n \in \N$, there exists a collection $\G$ of graphs on $n$ vertices, with
  \begin{equation}\label{eq:CT:triangles:size}
    |\G| \le n^{C n^{3/2}},
  \end{equation}
  such that
  \begin{itemize}
  \item[$(a)$] each $G \in \G$ contains fewer than $\eps n^3$ triangles;
  \item[$(b)$] each triangle-free graph on $n$ vertices is contained in some $G \in \G$. 
  \end{itemize}
\end{thm}

In order to motivate the statement of Theorem~\ref{thm:CT:triangles}, we will next present three simple applications: bounding the largest size of a triangle-free subgraph of the random graph $G(n,p)$, determining the typical structure of a (sparse) triangle-free graph, and proving that $G(n,p)$ cannot be partitioned into a bounded number of triangle-free graphs.

\subsection{Mantel's theorem in random graphs} \label{Mantel:sec}

The oldest result in extremal graph theory, which states that every graph on $n$ vertices with more than $n^2/4$ edges contains a triangle, was proved by Mantel~\cite{Ma07} in 1907. The corresponding problem in the random graph $G(n,p)$ was first studied by Frankl and R\"odl~\cite{FR}, who proved the following theorem (cf.~Theorem~\ref{thm:Turan:Gnp}).

\begin{thm}\label{thm:Mantel:Gnp}
  For every $\alpha > 0$, there exists $C > 0$ such that the following holds. If $p \ge C / \sqrt{n}$, then a.a.s.\ every subgraph $G \subset G(n,p)$ with 
  \[
    e(G) \ge \bigg( \frac{1}{2} + \alpha \bigg) p \binom{n}{2}
  \]
  contains a triangle. 
\end{thm}

As a simple first application of Theorem~\ref{thm:CT:triangles}, let us use it to prove Theorem~\ref{thm:Mantel:Gnp} under the marginally stronger assumption that $p \gg \log n / \sqrt{n}$. The proof exploits the following crucial property of $n$-vertex graphs with $o(n^3)$ triangles: each such graph has at most $\big(\frac{1}{2}+o(1)\big) \binom{n}{2}$ edges. This statement is made rigorous in the following supersaturation lemma for triangles, which can be proved by simply applying Mantel's theorem to each induced subgraph of $G$ with $O(1)$ vertices.

\begin{lemma}[Supersaturation for triangles]\label{lem:supersat:triangles}
  For every $\delta > 0$, there exists $\eps > 0$ such that the following holds. If $G$ is a graph on $n$ vertices with
  \[
    e(G) \ge \bigg( \frac{1}{4} + \delta \bigg) n^2,
  \]
  then $G$ has at least $\eps n^3$ triangles. 
\end{lemma} 

Applying Lemma~\ref{lem:supersat:triangles} with $\delta = \alpha/2$ and Theorem~\ref{thm:CT:triangles} with $\eps = \eps(\delta)$ given by the lemma, we obtain a family of containers $\G$ such that each $G \in \G$ has fewer than $\eps n^3$ triangles and thus
\[
  e(G) \le \bigg( \frac{1 + \alpha}{2} \bigg) \binom{n}{2}
\]
for every $G \in \G$. Since every triangle-free graph is a subgraph of some container, if $G(n,p)$ contains a triangle-free graph with $m$ edges, then in particular $e\big( G \cap G(n,p) \big) \ge m$ for some $G \in \G$. Noting that $e\big( G \cap G(n,p) \big) \sim \Bin\big( e(G), p \big)$, standard estimates on the tail of the binomial distribution yield
\[
  \Pr\bigg( e\big( G \cap G(n,p) \big) \ge \bigg( \frac{1}{2} + \alpha \bigg) p \binom{n}{2} \bigg) \le e^{- \beta pn^2},
\]
for some constant $\beta = \beta(\alpha) > 0$. Therefore, taking a union bound over all containers $G \in \G$ and using the bound~\eqref{eq:CT:triangles:size}, we have (using the notation of Theorem~\ref{thm:Turan:Gnp})
\begin{equation}
  \label{eq:Mantel:Gnp}
  \Pr\bigg( \ex\big( G(n,p), K_3 \big) \ge \bigg( \frac{1}{2} + \alpha \bigg) p {n \choose 2} \bigg) \le n^{O(n^{3/2})} \cdot e^{- \beta pn^2} \to 0
\end{equation}
as $n \to \infty$, provided that $p \gg \log n / \sqrt{n}$. This gives the conclusion of Theorem~\ref{thm:Mantel:Gnp} under a slightly stronger assumption on $p$. In Section~\ref{keylemma:sec}, we show how to remove the extra factor of $\log n$.

We remark here that Theorem~\ref{thm:Mantel:Gnp}, as well as numerous results of this type that now exist in the literature, cannot be proved using standard first moment estimates. Indeed, since there are at least $\binom{\lfloor n^2/4 \rfloor}{m}$ triangle-free graphs with $n$ vertices and $m$ edges, then letting $X_m$ denote the number of such graphs that are contained in $G(n,p)$, we have
\[
  \Ex[X_m] \ge p^m \binom{\lfloor n^2/4 \rfloor}{m} = \left(\frac{(e/2+o(1))p\binom{n}{2}}{m}\right)^m \gg 1
\]
if $m \le \big( e/2+o(1) \big) p\binom{n}{2} = o(n^2)$. This means that a first moment estimate would yield an upper bound on $\ex\big(G(n,p), K_3\big)$ that is worse than the trivial upper bound of $\big( 1+o(1) \big) p\binom{n}{2}$.


\subsection{The typical structure of a sparse triangle-free graph} 

A seminal theorem of Erd\H{o}s, Kleitman, and Rothschild~\cite{EKR} states that almost all triangle-free graphs are bipartite. Our second application of Theorem~\ref{thm:CT:triangles} is the following approximate version of this theorem for sparse graphs, first proved by {\L}uczak~\cite{Lu00}. Let us say that a graph $G$ is \emph{$t$-close to bipartite} if there exists a bipartite subgraph $G' \subset G$ with $e(G') \ge e(G) - t$.

\begin{thm}\label{thm:structure:trianglefree}
For every $\alpha > 0$, there exists $C > 0$ such that the following holds. If $m \ge C n^{3/2}$, then almost all triangle-free graphs with $n$ vertices and $m$ edges are $\alpha m$-close to bipartite.
\end{thm}

We will again (cf.~the previous subsection) prove this theorem under the marginally stronger assumption that $m \gg n^{3/2}\log n$. To do so, we will need a finer characterisation of graphs with $o(n^3)$ triangles that takes into account whether or not a graph is close to bipartite. Proving such a result is less straightforward than Lemma~\ref{lem:supersat:triangles}; for example, one natural proof combines the triangle removal lemma of Ruzsa and Szemer\'edi~\cite{RuSz} with the classical stability theorem of Erd\H{o}s and Simonovits~\cite{Er67, Si68}. However, an extremely simple, beautiful, and elementary proof was given recently by F\"uredi~\cite{Fur} (see also~\cite{BBCLMS}). 

\begin{lemma}[Robust stability for triangles]\label{lem:superstability:triangles}
For every $\delta > 0$, there exists $\eps > 0$ such that the following holds. If $G$ is a graph on $n$ vertices with
\[
e(G) \ge \bigg( \frac{1}{2} - \eps \bigg) {n \choose 2},
\]
then either $G$ is $\delta n^2$-close to bipartite or $G$ contains at least $\eps n^3$ triangles. 
\end{lemma} 


Applying Lemma~\ref{lem:superstability:triangles} with $\delta = \delta(\alpha) > 0$ sufficiently small and Theorem~\ref{thm:CT:triangles} with $\eps = \eps(\delta)$ given by the lemma, we obtain a family of containers $\G$ such that every $G \in \G$ is either $\delta n^2$-close to bipartite or
\begin{equation}\label{eq:container:with:few:edges}
e(G) \le \bigg( \frac{1}{2} - \eps \bigg) {n \choose 2}.
\end{equation}
Let us count those triangle-free graphs $H$ with $n$ vertices and $m$ edges that are not $\alpha m$-close to bipartite; note that each such graph is a subgraph of some container $G \in \G$. 

Suppose first that $G$ satisfies~\eqref{eq:container:with:few:edges}; in this case we  simply use the trivial bound 
\[
{e(G) \choose m} \le \binom{\left(\frac{1}{2}-\eps\right)\binom{n}{2}}{m} \le (1 - \eps)^m {n^2/4 \choose m}
\]
for the number of choices for $H \subset G$. On the other hand, if $G$ is $\delta n^2$-close to bipartite, then there is some bipartite $G' \subset G$ with $e(G') \ge e(G) - \delta n^2$. Since $e(H \cap G') \le (1-\alpha)m$ by our assumption on $H$, we bound the number of choices for $H$ by
\[
{e(G)-e(G') \choose \alpha m} {e(G) \choose (1 - \alpha) m} \le \binom{\delta n^2}{\alpha m} \binom{\binom{n}{2}}{(1-\alpha)m} \le 2^{-m} \binom{n^2/4}{m},
\]
provided that $\delta = \delta(\alpha)$ is sufficiently small. Summing over all choices of $G \in \G$ and using~\eqref{eq:CT:triangles:size}, it follows that if $m \gg n^{3/2} \log n$, then there are at most
\[
n^{O(n^{3/2})} \cdot (1 - \eps)^m \binom{n^2/4}{m} \ll \binom{\lfloor n^2/4 \rfloor}{m}
\]
triangle-free graphs $H$ with $n$ vertices and $m$ edges that are not $\alpha m$-close to bipartite. However, there are clearly at least ${\lfloor n^2/4 \rfloor \choose m}$ triangle-free graphs $H$ with $n$ vertices and $m$ edges, since every bipartite graph is triangle-free, so the conclusion of Theorem~\ref{thm:structure:trianglefree} holds when $m \gg n^{3/2} \log n$. We again postpone a discussion of how to remove the unwanted factor of $\log n$ to Section~\ref{keylemma:sec}. 

\subsection{Ramsey properties of sparse random graphs} 

A folklore fact that is presented in each introduction to Ramsey theory states that every $2$-colouring of the edges of $K_6$ contains a monochromatic triangle. With the aim of constructing a small $K_4$-free graph that has the same property, Frankl and R\"odl~\cite{FR} proved that if $p \gg 1/\sqrt{n}$, then a.a.s.\ every $2$-colouring of the edges of $G(n,p)$ contains a monochromatic triangle. Ramsey properties of random graphs were later thorougly investigated by R\"odl and Ruci\'nski~\cite{RoRu93, RoRu94, RoRu95}. The following theorem is the main result of~\cite{RoRu94}.

\begin{thm}
  \label{thm:triangle-Ramsey}
  For every $r \in \N$, there exists $C > 0$ such that the following holds. If $p \gg C/\sqrt{n}$, then a.a.s.\ every $r$-colouring of the edges of $G(n,p)$ contains a monochromatic triangle.
\end{thm}

We will present a simple proof of this theorem that was discovered recently by Nenadov and Steger~\cite{NS}. For the sake of simplicity, we will again use the marginally stronger assumption that $p \gg \log n / \sqrt{n}$. The proof exploits the following property of $n$-vertex graphs with $o(n^3)$ triangles: the union of any bounded number of such graphs cannot cover a $\big( 1 - o(1) \big)$-proportion of the edges of $K_n$. This property is a straightforward corollary of the following lemma, which can be proved by applying Ramsey's theorem to the colourings induced by all subsets of $V(K_n)$ of size $O(1)$.

\begin{lemma}\label{lemma:Ramsey-supersat}
  For every $r \in \N$, there exist $n_0$ and $\eps > 0$ such that for all $n \ge n_0$, every $(r+1)$-colouring of the edges of $K_n$ contains at least $(r+1) \eps n^3$ monochromatic triangles.
\end{lemma}

Applying Theorem~\ref{thm:CT:triangles} with $\eps = \eps(r)$ given by the lemma, we obtain a family of containers $\G$ such that every $G \in \G$ has fewer than $\eps n^3$ triangles. If $G(n,p)$ does not have the desired Ramsey property, then there are triangle-free graphs $H_1, \dotsc, H_r$ such that $H_1 \cup \dotsc \cup H_r = G(n,p)$. It follows that $G(n,p) \subset G_1 \cup \dotsc \cup G_r$, where each $G_i \in \G$ is a container for $H_i$. Since each $G_i$ has fewer than $\eps n^3$ triangles, then Lemma~\ref{lemma:Ramsey-supersat} implies that $K_n \setminus (G_1 \cup \dotsc \cup G_r)$ contains at least $\eps n^3$ triangles.\footnote{To see this, consider an $(r+1)$-colouring of the edges of $K_n$ that assigns to each edge $e \in G_1 \cup \dotsc \cup G_r$ some colour~$i$ such that $e \in G_i$ and assigns colour $r+1$ to all edges of $K_n \setminus (G_1 \cup \dotsc \cup G_r)$.} Since each edge of $K_n$ belongs to fewer than $n$ triangles, we must have $e\big(K_n \setminus (G_1 \cup \dotsc \cup G_r)\big) \ge \eps n^2$. Consequently, for each fixed $G_1, \dotsc, G_r \in \G$,
\[
\Pr\big( G(n,p) \subset G_1 \cup \dotsc \cup G_r \big) = (1-p)^{e(K_n \setminus (G_1 \cup \cdots \cup G_r))} \le (1-p)^{\eps n^2} \le e^{- \eps pn^2}.
\]
Taking a union bound over all $r$-tuples of containers, we conclude that
\[
\Pr\big( \text{$G(n,p)$ admits a `bad' $r$-colouring} \big) \le n^{O(n^{3/2})} \cdot e^{- \eps pn^2} \to 0
\]
as $n \to \infty$, provided that $p \gg \log n / \sqrt{n}$. As before, the unwanted factor of $\log n$ can be removed with a somewhat more careful analysis that we shall discuss in Section~\ref{keylemma:sec}.

\subsection{An application in discrete geometry}\label{geometry:sec}

In order to give some idea of the flexibility of the container method, we will next present a somewhat more elaborate application of the container lemma for 3-uniform hypergraphs, which was discovered recently by Balogh and Solymosi~\cite{BS}, to the following question posed by Erd\H{o}s~\cite{Er88}. Given $n$ points in the Euclidean plane $\R^2$, with at most three on any line, how large a subset are we guaranteed to find in general position (i.e., with at most two on any line)? F\"uredi~\cite{Fur91} proved that one can always find such a subset of size $\Omega\big(\sqrt{n\log n}\big)$ and gave a construction (which relied on the density Hales--Jewett theorem of Furstenberg and Katznelson~\cite{FuKa}) in which the largest such set has size $o(n)$. Using the method of hypergraph containers, Balogh and Solymosi~\cite{BS} obtained the following stronger upper bound. 

\begin{thm}\label{thm:geometry}
There exists a set $S \subset \R^2$ of size $n$, containing no four points on a line, such that every subset of $S$ of size $n^{5/6+o(1)}$ contains three points on a line.
\end{thm}

The key idea in~\cite{BS} is to first construct a set $P$ of points that contains `few' collinear quadruples, but such that every `large' subset of $P$ contains `many' collinear triples. Then a random subset $R$ of $P$ of a carefully chosen density will typically contain only $o(|R|)$ collinear quadruples, since the density is not too large and there are few collinear quadruples. On the other hand, every subset of $R$ with more than $|R|^{5/6+o(1)}$ elements will still contain a collinear triple; this follows from the hypergraph container lemma, as large sets contain many collinear triples and the density is not too small. Removing one element from each collinear quadruple in $R$ gives the desired set~$A$.


Formally, we first define the following 3-uniform hypergraph $\cH$. We let $V(\cH) = [m]^3$ (so the vertices are lattice points in $\R^3$) and let $E(\cH)$ be the collection of triples of points that lie on a common line. Thus, a subset of $V(\cH)$ is in general position if and only if it is an independent set of $\cH$. The following lemma was proved in~\cite{BS}. 

\begin{lemma}[Supersaturation for collinear triples]\label{lem:supersat:triples:on:lines}
  For every $0 < \gamma < 1/2$ and every $S \subset [m]^3$ of size at least $m^{3-\gamma}$, there exist at least $m^{6 - 4\gamma - o(1)}$ collinear triples of points in $S$. 
\end{lemma} 

We now repeatedly apply the hypergraph container lemma for 3-uniform hypergraphs to subhypergraphs of $\cH$. Suppose that $s \ge m^{8/3+o(1)}$ and let $S \subset [m]^3$ be an arbitrary $s$-element set. Lemma~\ref{lem:supersat:triples:on:lines} gives
\[
  e\big( \cH[S] \big) \ge s^4/m^{6+o(1)} \qquad \text{and} \qquad \Delta_2\big( \cH[S] \big) \le \Delta_2(\cH) \le m.
\]
Moreover, it is not difficult to deduce that there exists a subhypergraph $\cH' \subset \cH[S]$ with
\[
  v(\cH') = |S| = s, \quad e(\cH') = s^4 / m^{6+o(1)}, \quad \text{and} \quad \Delta_1(\cH') = O\big(e(\cH') / v(\cH')\big).
\]
We may therefore apply the container lemma for 3-uniform hypergraphs to $\cH'$ to obtain a collection $\C$ of at most $\exp\big( m^{3+o(1)}/\sqrt{s} \big)$ subsets of $S$ with the following properties:
\begin{enumerate}
\item[$(a)$] Every set of points of $S$ in general position is contained in some $C \in \C$,\smallskip
\item[$(b)$] Each $C \in \C$ has size at most $(1 - \delta) |S|$.
\end{enumerate}
Starting with $S = [m]^3$ and iterating this process for $O(\log m)$ steps, we obtain the following container theorem for sets of points in general position.

\begin{thm}\label{thm:CT:points:gen:position}
  For each $m \in \N$, there exists a collection $\C$ of subsets of $[m]^3$ with
  \begin{equation}\label{eq:CT:points:gen:position:size}
    |\C| \le \exp\big( m^{5/3 + o(1)} \big)
  \end{equation}
  such that
  \begin{itemize}
  \item[$(a)$] $|C| \le m^{8/3+o(1)}$ for each $C \in \C$;
  \item[$(b)$] each set of points of $[m]^3$ in general position is contained in some $C \in \C$. 
  \end{itemize}
\end{thm}

Now, let $p = m^{-1+o(1)}$ and consider a $p$-random subset $R \subset [m]^3$, that is, each element of $[m]^3$ is included in $R$ independently at random with probability $p$. Since $[m]^3$ contains $m^{6+o(1)}$ sets of four collinear points\footnote{This is because there are $O(m^6/t^4)$ lines in $\R^3$ that contain more than $t$ points of $[m]^3$.}, it follows that, with high probability, $|R| = pm^{3+o(1)} = m^{2+o(1)}$ and $R$ contains $p^4 m^{6+o(1)} = o(|R|)$ collinear 4-tuples. Moreover, since $|C| \le m^{8/3+o(1)}$ for each $C \in \C$, it follows from~\eqref{eq:CT:points:gen:position:size} and standard estimates on the tail of the binomial distribution that with high probability we have $|R \cap C| \le m^{5/3 + o(1)}$ for every $C \in \C$. In particular, removing one element from each collinear 4-tuple in $R$ yields a set $A \subset [m]^3$ of size $m^{2+o(1)}$ with no collinear 4-tuple and containing no set of points in general position of size larger than $m^{5/3 + o(1)}$. Finally, project the points of~$A$ to the plane in such a way that collinear triples remain collinear, and no new collinear triple is created. In this way, we obtain a set of $n = m^{2+o(1)}$ points in the plane, no four of them on a line, such that no set of size greater than $n^{5/6+o(1)} = m^{5/3+o(1)}$ is in general position, as required.



\section{The key container lemma}\label{keylemma:sec}

In this section, we state a container lemma for hypergraphs of arbitrary uniformity. The version of the lemma stated below, which comes from~\cite{MSS}, differs from the statement originally proved by the authors of this survey~\cite[Proposition~3.1]{BMS} only in that the dependencies between the various constants have been made more explicit here; a careful analysis of the proof of~\cite[Proposition~3.1]{BMS} will yield this slightly sharper statement.\footnote{A complete proof of the version of the container lemma stated here can be found in~\cite{MSS}.} Let us recall that for a hypergraph $\cH$ and an integer $\ell$, we write $\Delta_\ell(\cH)$ for the maximum degree of a set of $\ell$ vertices of $\cH$, that is,
\[
  \Delta_\ell(\cH) = \max\big\{ d_\cH(A) \scolon A \subset V(\cH), \, |A| = \ell \big\},
\]
where $d_\cH(A) = \big| \big\{ B \in E(\cH) \scolon A \subset B \big\} \big|$, and $\I(\cH)$ for the collection of independent sets of $\cH$. The lemma states, roughly speaking, that each independent set $I$ in a uniform hypergraph $\cH$ can be assigned a \emph{fingerprint} $S \subset I$ in such a way that all sets with the same fingerprint are contained in a single set $C = f(S)$, called a \emph{container}, whose size is bounded away from $v(\cH)$. More importantly, the sizes of these fingerprints (and hence also the number of containers) can be bounded from above (in an optimal way!) by basic parameters of $\cH$.

\begin{HCL}
  Let $k \in \N$ and set $\delta = 2^{-k(k+1)}$. Let $\cH$ be a $k$-uniform hypergraph and suppose that
  \begin{equation}\label{eq:containers:condition}
    \Delta_\ell(\cH) \le \bigg( \frac{b}{v(\cH)} \bigg)^{\ell-1} \, \frac{e(\cH)}{r}
  \end{equation}
  for some $b,r \in \N$ and every $\ell \in \{1, \dotsc, k\}$. Then there exists a collection $\C$ of subsets of $V(\cH)$ and a function $f \colon \P\big( V(\cH) \big) \to \C$ such that:
  \begin{enumerate}
  \item[$(a)$] for every $I \in \I(\cH)$, there exists $S \subset I$ with $|S| \le (k-1)b$ and $I \subset f(S)$;
  \item[$(b)$] $|C| \le v(\cH) - \delta r$ for every $C \in \C$. 
  \end{enumerate}
\end{HCL}

The original statement of the container lemma~\cite[Proposition~3.1]{BMS} had $r = v(\cH) / c$ for some constant $c$, since this choice of parameters is required in most standard applications. In particular, the simple container lemma for $3$-uniform hypergraphs presented in Section~\ref{basic:sec} is easily derived from the above statement by letting $b = v(\cH) / (2\sqrt{d})$ and $r = v(\cH) / (6c)$, where $d = 3e(\cH)/v(\cH)$ is the average degree of $\cH$. There are, however, arguments that benefit from setting $r = o(v(\cH))$; we present one of them in Section~\ref{sec:many-colours}.

Even though the property $|C| \le v(\cH) - \delta r$ that is guaranteed for all containers $C \in \C$ seems rather weak at first sight, it can be easily strengthened with repeated applications of the lemma. In particular, if for some hypergraph $\cH$, condition~\eqref{eq:containers:condition} holds (for all $\ell$) with some $b = o(v(\cH))$ and $r = \Omega(v(\cH))$, then recursively applying the lemma to subhypergraphs of $\cH$ induced by all the containers $C$ for which $e(\cH[C]) \ge \eps e(\cH)$ eventually produces a collection $\C$ of containers indexed by sets of size $O(b)$ such that $e(\cH[C]) < \eps e(\cH)$ for every $C \in \C$. This is precisely how (in Section~\ref{basic:sec}) we derived Theorem~\ref{thm:CT:triangles} from the container lemma for $3$-uniform hypergraphs. For a formal argument showing how such a family of `tight' containers may be constructed, we refer the reader to~\cite{BMS}.

One may thus informally say that the hypergraph container lemma provides a covering of the family of all independent sets of a uniform hypergraph with `few' sets that are `almost independent'. In many natural settings, these almost independent sets closely resemble truly independent sets. In some cases, this is a straightforward consequence of corresponding removal lemmas. A more fundamental reason is that many sequences of hypergraphs $\cH_n$ of interest possess the following self-similarity property: For all (or many) pairs $m$ and $n$ with $m < n$, the hypergraph $\cH_n$ admits a very uniform covering by copies of $\cH_m$. For example, this is the case when $\cH_n$ is the hypergraph encoding triangles in $K_n$, simply because every $m$-element set of vertices of $K_n$ induces $K_m$. Such self-similarity enables one to use elementary averaging arguments to characterise almost independent sets; for example, the standard proof of Lemma~\ref{lem:supersat:triangles} uses such an argument.

The fact that the fingerprint $S$ of each independent set $I \in \I(\cH)$ is a subset of $I$ is not merely a by-product of the proof of the hypergraph container lemma. On the contrary, it is an important property of the family of containers that can be often exploited to make union bound arguments tighter. This is because each $I \in \I(\cH)$ is sandwiched between $S$ and $f(S)$ and consequently when enumerating independent sets one may use a union bound over all fingerprints $S$ and enumerate only over the sets $I \setminus S$ (which are contained in $f(S)$). In particular, such finer arguments can be used to remove the superfluous logarithmic factor from the assumptions of the proofs outlined in Section~\ref{basic:sec}. 
For example, in the proof of Theorem~\ref{thm:Mantel:Gnp} presented in Section~\ref{Mantel:sec}, the fingerprints of triangle-free subgraphs of $K_n$ form a family $\S$ of $n$-vertex graphs, each with at most $C_\eps n^{3/2}$ edges. Setting $m = \big( \frac{1}{2} + \alpha \big) p \binom{n}{2}$, this allows us to replace~\eqref{eq:Mantel:Gnp} with the following estimate:
\begin{equation}\label{eq:fingerprint}
\Pr\Big( \ex\big( G(n,p), K_3 \big) \ge m \Big) \, \le \, \sum_{S \in \S} \Pr\Big( S \subset G(n,p) \textup{ and } e\big( \big(f(S)\setminus S\big) \cap G(n,p) \big) \ge m - |S| \Big).
\end{equation}
Since the two events in the right-hand side of~\eqref{eq:fingerprint} concern the intersections of $G(n,p)$ with two disjoint sets of edges of $K_n$, they are independent. If $p \gg n^{-1/2}$, then $|S| \ll p \binom{n}{2}$ and consequently, recalling that $e(f(S)) \le \big(\frac{1+\alpha}{2}\big) \binom{n}{2}$, we may bound the right-hand side of~\eqref{eq:fingerprint} from above by
\[
  \sum_{S \in \S} p^{|S|} e^{-\beta pn^2} \le \sum_{s \le C_\eps n^{3/2}} \binom{\binom{n}{2}}{s} \cdot p^s e^{-\beta pn^2} \le \sum_{s \le C_\eps n^{3/2}} \left(\frac{e\binom{n}{2}p}{s}\right)^s e^{-\beta pn^2} \le e^{-\beta pn^2/2}
\]
for some $\beta = \beta(\alpha) > 0$.

Finally, what is the intuition behind condition~\eqref{eq:containers:condition}? A natural way to define $f(S)$ for a given (independent) set $S$ is to let $f(S) = V(\cH) \setminus X(S)$, where $X(S)$ comprises all vertices $v$ such that $A \subset S \cup \{v\}$ for some $A \in E(\cH)$. Indeed, every independent set $I$ that contains $S$ must be disjoint from $X(S)$. (In reality, the definition of $X(S)$ is -- and has to be -- more complicated than this, and some vertices are placed in $X(S)$ simply because they do not belong to $S$.) Suppose, for the sake of argument, that $S$ is a random set of $b$ vertices of $\cH$. Letting $\tau = b/v(\cH)$, we have
\begin{equation}
  \label{eq:fingerprint-intuition}
  \Ex\big[|X(S)|\big] \le \sum_{A \in E(\cH)} \Pr\big(|A \cap S| = k-1\big) \le k \cdot \tau^{k-1} \cdot e(\cH).
\end{equation}
Since we want $X(S)$ to have at least $\delta r$ elements for every fingerprint $S$, it seems reasonable to require that
\[
  \Delta_k(\cH) = 1 \le \frac{k}{\delta} \cdot \tau^{k-1} \cdot \frac{e(\cH)}{r},
\]
which is, up to a constant factor, condition~\eqref{eq:containers:condition} with $\ell = k$. For some hypergraphs $\cH$ however, the first inequality in~\eqref{eq:fingerprint-intuition} can be very wasteful, since some $v \in X(S)$ may have many $A \in E(\cH)$ such that $A \subset S \cup \{v\}$.  This can happen if for some $\ell \in \{1, \dotsc, k-1\}$, there is an $\ell$-uniform hypergraph $\G$ such that each edge of $\cH$ contains an edge of $\G$; note that $e(\G)$ can be as small as $e(\cH) / \Delta_\ell(\cH)$. Our assumption implies that $\I(\G) \subset \I(\cH)$ and thus, letting $Y(S)$ be the set of all vertices $w$ such that $B \subset S \cup \{w\}$ for some $B \in E(\G)$, we have $X(S) \subset Y(S)$. In particular, we want $Y(S)$ to have at least $\delta r$ elements for every fingerprint $S$ of an independent set $I \in \I(\G)$. Repeating~\eqref{eq:fingerprint-intuition} with $X$ replaced by $Y$, $\cH$ replaced by $\G$, and $k$ replaced by $\ell$, we arrive at the inequality
\[
  \delta r \le \ell \cdot \tau^{\ell-1} \cdot e(\G) = \ell \cdot \tau^{\ell-1} \cdot \frac{e(\cH)}{\Delta_\ell(\cH)},
\]
which is, up to a constant factor, condition~\eqref{eq:containers:condition}. 

One may further develop the above argument to show that condition~\eqref{eq:containers:condition} is asymptotically optimal, at least when $r = \Omega(v(\cH))$. Roughly speaking, one can construct $k$-uniform hypergraphs that have $\binom{(1-o(1)) v(\cH)}{m}$ independent $m$-sets for every $m = o(b)$, where $b$ is minimal so that condition~\eqref{eq:containers:condition} holds, whereas the existence of containers of size at most $(1-\delta)  v(\cH)$ indexed by fingerprints of size $o(b)$ would imply that the number of such sets is at most $\binom{(1-\eps) v(\cH)}{m}$ for some constant $\eps > 0$.


\section{Counting $H$-free graphs}\label{counting:Hfree:sec}

How many graphs are there on $n$ vertices that do not contain a copy of $H$? An obvious lower bound is $2^{\ex(n,H)}$, since each subgraph of an $H$-free graph is also $H$-free. For non-bipartite graphs, this is not far from the truth.  Writing $\cF_n(H)$ for the family of $H$-free graphs on $n$ vertices, if $\chi(H) \ge 3$, then
\begin{equation}\label{eq:EFR}
|\cF_n(H)| = 2^{(1 + o(1))\ex(n,H)}
\end{equation}
as $n \to \infty$, as was first shown by Erd\H{o}s, Kleitman, and Rothschild~\cite{EKR} (when $H$ is a complete graph) and then by Erd\H{o}s, Frankl, and R\"odl~\cite{EFR}. For bipartite graphs, on the other hand, the problem is much more difficult. In particular, the following conjecture (first stated in print in~\cite{KW82}), which played a major role in the development of the container method, remains open.

\begin{conj}\label{conj:counting}
  For every bipartite graph $H$ that contains a cycle, there exists $C > 0$ such that
  \[
    |\cF_n(H)| \le 2^{C\ex(n,H)}
  \]
  for every $n \in \N$.
\end{conj}

The first significant progress on Conjecture~\ref{conj:counting} was made by Kleitman and Winston~\cite{KW82}. Their proof of the case $H = C_4$ of the conjecture introduced (implicitly) the container method for graphs. Nevertheless, it took almost thirty years\footnote{An unpublished manuscript of Kleitman and Wilson from 1996 proves that $|\cF_n(C_6)| = 2^{O(\ex(n,C_6))}$.} until their theorem was generalized to the case $H = K_{s,t}$, by Balogh and Samotij~\cite{BSmm,BSst}, and then (a few years later) to the case $H = C_{2k}$, by Morris and Saxton~\cite{MS}. More precisely, it was proved in~\cite{BSst,MS} that
\[
  |\cF_n(K_{s,t})| = 2^{O(n^{2-1/s})} \qquad \textup{and} \qquad |\cF_n(C_{2k})| = 2^{O(n^{1+1/k})}
\]
for every $2 \le s \le t$ and every $k \ge 2$, which implies Conjecture~\ref{conj:counting} when $t > (s-1)!$ and $k \in \{2,3,5\}$, since in these cases it is known that $\ex(n,K_{s,t}) = \Theta(n^{2-1/s})$ and $\ex(n,C_{2k}) = \Theta(n^{1+1/k})$. 

Very recently, Ferber, McKinley, and Samotij~\cite{FMS}, inspired by a similar result of Balogh, Liu, and Sharifzadeh~\cite{BLS} on sets of integers with no $k$-term arithmetic progression, found a very simple proof of the following much more general theorem.

\begin{thm}\label{thm:FMS}
  Suppose that $H$ contains a cycle. If $\ex(n,H) = O(n^\alpha)$ for some constant $\alpha$, then 
  \[
    |\cF_n(H)| = 2^{O(n^\alpha)}.
  \]
\end{thm}

Note that Theorem~\ref{thm:FMS} resolves Conjecture~\ref{conj:counting} for every $H$ such that $\ex(n,H) = \Theta(n^\alpha)$ for some constant $\alpha$. Moreover, it was shown in~\cite{FMS} that the weaker assumption that $\ex(n,H) \gg n^{2 - 1/m_2(H) + \eps}$ for some $\eps > 0$ already implies that the assertion of Conjecture~\ref{conj:counting} holds for infinitely many $n$; we refer the interested reader to~\cite{FMS} for details. Let us also note here that, while it is natural to suspect that in fact the stronger bound~\eqref{eq:EFR} holds for all graphs $H$ that contain a cycle, this is false for $H = C_6$, as was shown by Morris and Saxton~\cite{MS}. However, it may still hold for $H = C_4$ and it would be very interesting to determine whether or not this is indeed the case.

The proof of Theorem~\ref{thm:FMS} for general $H$ is somewhat technical, so let us instead sketch the proof in the case $H = C_4$. In this case, the proof combines the hypergraph container lemma stated in the previous section with the following supersaturation lemma.

\begin{lemma}\label{lem:C4:super}
  There exist constants $\beta > 0$ and $k_0 \in \N$ such that the following holds for every $k \ge k_0$ and every $n \in \N$. Given a graph~$G$ with~$n$ vertices and $k \cdot \ex(n,C_4)$ edges, there exists a collection $\cH$ of at least $\beta k^5 \cdot \ex(n,C_4)$ copies of $C_4$ in $G$ that satisfies:
  \begin{itemize}
  \item[$(a)$] Each edge belongs to at most $k^4$ members of $\cH$.
  \item[$(b)$] Each pair of edges is contained in at most $k^2$ members of $\cH$.
  \end{itemize}
\end{lemma}

The proof of Lemma~\ref{lem:C4:super} employs several simple but important ideas that can be used in a variety of other settings, so let us sketch the details. The first key idea, which was first used in~\cite{MS}, is to build the required family $\cH$ one $C_4$ at a time. Let us say that a collection~$\cH$ of copies of $C_4$ is \emph{legal} if it satisfies conditions $(a)$ and $(b)$ and suppose that we have already found a legal collection~$\cH_m$ of $m$ copies of $C_4$ in $G$. Note that we are done if $m \ge \beta k^5 \cdot \ex(n,C_4)$, so let us assume that the reverse inequality holds and construct a legal collection~$\cH_{m+1} \supset \cH_m$ of $m+1$ copies of $C_4$ in $G$.

We claim that there exists a collection $\A_m$ of $\beta k^5 \cdot \ex(n,C_4)$ copies of $C_4$ in $G$, any of which can be added to $\cH_m$ without violating conditions $(a)$ and $(b)$, that is, such that $\cH_m \cup \{ C \}$ is legal for any $C \in \A_m$. (Let us call these \emph{good} copies of $C_4$.) Since $m < \beta k^5 \cdot \ex(n,C_4)$, then at least one element of $\A_m$ is not already in $\cH_m$, so this will be sufficient to prove the lemma.

To find $\A_m$, observe first that (by simple double-counting) at most $4\beta k \cdot \ex(n,C_4)$ edges of $G$ lie in exactly $k^4$ members of $\cH_m$ and similarly at most $6\beta k^3 \cdot \ex(n,C_4)$ pairs of edges of $G$ lie in exactly $k^2$ members of $\cH_m$.  Now, consider a random subset $A \subset V(G)$ of size $pn$, where $p = D/k^2$ for some large constant $D$. Typically $G[A]$ contains about $p^2 k \cdot \ex(n,C_4)$ edges. After removing from $G[A]$ all \emph{saturated} edges (i.e., those belonging to $k^4$ members of $\cH_m$) and one edge from each \emph{saturated} pair (i.e., pair of edges that is contained in $k^2$ members of $\cH_m$), we expect to end up with at least
\begin{equation*}
p^2 k \cdot \ex(n,C_4) - 4\beta p^2 k \cdot \ex(n,C_4) - 6\beta p^3 k^3 \cdot \ex(n,C_4) \ge \frac{p^2 k \cdot \ex(n,C_4)}{2} \ge 2 \cdot \ex( pn, C_4)
\end{equation*}
edges, where the first inequality follows since $p = D / k^2$ and $\beta$ is sufficiently small, and the second holds because $\ex(n,C_4) = \Theta(n^{3/2})$ and $D$ is sufficiently large.  Finally, observe that any graph on $pn$ vertices with at least $2 \cdot \ex( pn, C_4)$ edges contains at least 
$$\ex\big( pn, C_4 \big) \, = \, \Omega\Big( p^{3/2} \cdot \ex(n,C_4) \Big)$$ 
copies of~$C_4$. But each copy of $C_4$ in $G$ was included in the random subgraph $G[A]$ with probability at most $p^4$ and hence (with a little care) one can show that there must exist at least $\Omega\big( p^{-5/2} \cdot \ex(n, C_4) \big)$ copies of $C_4$ in $G$ that avoid all saturated edges and pairs of edges. Since $p^{-5/2} = k^5/D^{5/2}$ and $\beta$ is sufficiently small, we have found $\beta k^5 \cdot \ex(n, C_4)$ good copies of $C_4$ in $G$, as required.

We now show how  one may combine Lemma~\ref{lem:C4:super} and the hypergraph container lemma to construct families of containers for $C_4$-free graphs. Let $\beta$ and $k_0$ be the constants from the statement of Lemma~\ref{lem:C4:super} and assume that $G$ is an $n$-vertex graph with at least $k \cdot \ex(n, C_4)$ and at most $2k \cdot \ex(n, C_4)$ edges, where $k \ge k_0$. Denote by $\cH_G$ the 4-uniform hypergraph with vertex set $E(G)$, whose edges are the copies of $C_4$ in $G$ given by Lemma~\ref{lem:C4:super}. Since
\[
  v(\cH_G) = e(G), \quad e(\cH_G) \ge \beta k^5 \cdot \ex(n,C_4), \quad \Delta_1(\cH_G) \le k^4, \quad \Delta_2(\cH_G) \le k^2,
\]
and $\Delta_3(\cH_G) = \Delta_4(\cH_G) = 1$, the hypergraph $\cH_G$ satisfies the assumptions of the container lemma with $r = \beta k \cdot \ex(n, C_4)$ and $b = 2k^{-1/3} \cdot \ex(n, C_4)$. Consequently, there exist an absolute constant~$\delta$ and a collection $\C$ of subgraphs of $G$ with the following properties:
\begin{enumerate}
\item[$(a)$] every $C_4$-free subgraph of $G$ is contained in some $C \in \C$,
\item[$(b)$] each $C \in \C$ has at most $(1 - \delta) e(G)$ edges,
\end{enumerate}
and moreover
\[
  |\C| \le \sum_{s=0}^{3b} \binom{e(G)}{s} \le \left(\frac{e(G)}{b}\right)^{3b} \le k^{4b} \le \exp\Big( 8 k^{-1/3} \log k \cdot \ex(n,C_4) \Big).
\]
Note that we have just replaced a single container for the family of $C_4$-free subgraphs of $G$ (namely $G$ itself) with a small collection of containers for this family, each of which is somewhat smaller than~$G$. Since every $C_4$-free graph with $n$ vertices is contained in $K_n$, by repeatedly applying this `breaking down' process, we obtain the following container theorem for $C_4$-free graphs.

\begin{thm}\label{thm:C4containers:weak}
  There exist constants $k_0 > 0$ and $C > 0$ such that the following holds for all $n \in \N$ and $k \ge k_0$. There exists a collection~$\G(n,k)$ of at most
  \[
    \exp\bigg( \frac{C \log k}{k^{1/3}} \cdot \ex(n,C_4) \bigg)
  \]
  graphs on $n$ vertices such that
  \[
    e(G) \le k \cdot \ex(n,C_4)
  \]
  for every $G \in \G(n,k)$ and every $C_4$-free graph on $n$ vertices is a subgraph of some~$G \in \G(n,k)$. 
\end{thm}

To obtain the claimed upper bound on $|\G(n,k)|$, note that if $k \cdot \ex(n,C_4) \ge \binom{n}{2}$ then we may take $\G(n,k) = \{K_n\}$, and otherwise the argument presented above yields
\[
  |\G(n,k)| \le \big|\G\big(n, k/(1-\delta)\big)\big| \cdot \exp\Big(8k^{-1/3}\log k \cdot \ex(n, C_4)\Big).
\]
In particular, applying Theorem~\ref{thm:C4containers:weak} with $k = k_0$, we obtain a collection of $2^{O(\ex(n,C_4))}$ containers for $C_4$-free graphs on $n$ vertices, each with $O\big( \ex(n,C_4) \big)$ edges. This immediately implies that Conjecture~\ref{conj:counting} holds for $H = C_4$. The proof for a general graph $H$ (under the assumption that $\ex(n,H) = \Theta(n^\alpha)$ for some $\alpha \in (1,2)$) is similar, though the details are rather technical. 

\subsection{Tur\'an's problem in random graphs}

Given that the problem of estimating $|\cF_n(H)|$ for bipartite graphs $H$ is notoriously difficult, it should not come as a surprise that determining the typical value of the Tur\'an number $\ex\big(G(n,p), C_4\big)$ for bipartite $H$ also poses considerable challenges. Compared to the non-bipartite case, which was essentially solved by Conlon--Gowers~\cite{CG} and Schacht~\cite{Sch}, see Theorem~\ref{thm:Turan:Gnp}, the typical behaviour of $\ex\big(G(n,p), H \big)$ for bipartite graphs $H$ is much more subtle.

For simplicity, let us restrict our attention to the case $H = C_4$. Recall from Theorem~\ref{thm:Turan:Gnp} that the typical value of $\ex\big(G(n,p), C_4\big)$ changes from $\big( 1 + o(1) \big) p\binom{n}{2}$ to $o(pn^2)$ when $p = \Theta(n^{-2/3})$, as was first proved by Haxell, Kohayakawa, and {\L}uczak~\cite{HKL-even}. However, already several years earlier F\"uredi~\cite{Fur91} used the method of Kleitman and Winston~\cite{KW82} to prove\footnote{To be precise, F\"uredi proved that, if $m \ge 2n^{4/3} (\log n)^2$, then there are at most $(4n^3 / m^2)^m$ $C_4$-free graphs with $n$ vertices and $m$ edges, which implies the upper bounds in Theorem~\ref{thm:randomturan}. For the lower bounds, see~\cite{KKS,MS}.} the following much finer estimates of this extremal number for $p$ somewhat above the threshold.

\begin{thm}\label{thm:randomturan}
  Asymptotically almost surely,
  \[
    \ex \big( G(n,p), C_4 \big) = 
    \begin{cases}
      \big( 1 + o(1) \big) p {n \choose 2} & \textup{if $n^{-1} \ll p \ll n^{-2/3}$}, \\
      n^{4/3} (\log n)^{O(1)} & \textup{if $n^{-2/3} \le p \le n^{-1/3} (\log n)^4$}, \\
      \Theta\big( \sqrt{p} \cdot n^{3/2} \big) & \textup{if $p \ge n^{-1/3} (\log n)^4$}.
    \end{cases}
  \]
\end{thm}

We would like to draw the reader's attention to the (somewhat surprising) fact that in the middle range $n^{-2/3+o(1)} \le p \le n^{-1/3+o(1)}$, the typical value of $\ex \big( G(n,p), C_4 \big)$ stays essentially constant. A similar phenomenon has been observed in random Tur\'an problems for other forbidden bipartite graphs (even cycles~\cite{KKS, MS} and complete bipartite graphs~\cite{MS}) as well as Tur\'an-type problems in additive combinatorics~\cite{DeKoLeRoSa-Bh, DeKoLeRoSa-B3}. It would be very interesting to determine whether or not a similar `long flat segment' appears in the graph of $p \mapsto \ex\big( G(n,p), H \big)$ for every bipartite graph $H$. We remark that the lower bound in the middle range is given (very roughly speaking) by taking a random subgraph of $G(n,p)$ with density $n^{-2/3+o(1)}$ and then finding\footnote{One easy way to do this is simply to remove one edge from each copy of $C_4$. A more efficient method, used by Kohayakawa, Kreuter, and Steger~\cite{KKS} to improve the lower bound by a polylogarithmic factor, utilizes a version of the general result of~\cite{AKPSS} on independent sets in hypergraphs obtained in~\cite{DuLeRo95}; see also~\cite{FMS}.} a large $C_4$-free subgraph of this random graph; the lower bound in the top range is given by intersecting $G(n,p)$ with a suitable blow-up of an extremal $C_4$-free graph and destroying any $C_4$s that occur; see~\cite{KKS, MS} for details.

Even though Theorem~\ref{thm:C4containers:weak} immediately implies that $\ex\big( G(n,p), C_4\big) = o(pn^2)$ if $p \gg n^{-2/3} \log n$, it is not strong enough to prove Theorem~\ref{thm:randomturan}. A stronger container theorem for $C_{2\ell}$-free graphs (based on a supersaturation lemma that is sharper than Lemma~\ref{lem:C4:super}) was obtained in~\cite{MS}. In the case $\ell = 2$, the statement is as follows.

\begin{thm}\label{thm:C4containers:turan}
  There exist constants $k_0 > 0$ and $C > 0$ such that the following holds for all $n \in \N$ and $k_0 \le k \le n^{1/6} / \log n$. There exists a collection~$\G(n,k)$ of at most
  \[
    \exp\bigg( \frac{C \log k}{k} \cdot \ex(n,C_4) \bigg)
  \]
  graphs on $n$ vertices such that
  \[
    e(G) \le k \cdot \ex(n,C_4)
  \]
  for every $G \in \G(n,k)$ and every $C_4$-free graph on $n$ vertices is a subgraph of some~$G \in \G(n,k)$. 
\end{thm}

Choosing $k$ to be a suitable function of $p$, it is straightforward to use Theorem~\ref{thm:C4containers:turan} to prove a slightly weaker version of Theorem~\ref{thm:randomturan}, with an extra factor of $\log n$ in the upper bound on $\ex \big( G(n,p), C_4 \big)$. As usual, this logarithmic factor can be removed via a more careful application of the container method, using the fact that the fingerprint of an independent set is contained in it, cf.\ the discussion in Section~\ref{keylemma:sec}; see~\cite{MS} for the details. However, we are not able to determine the correct power of $\log n$ in $\ex \big( G(n,p), C_4 \big)$ in the middle range $n^{-2/3 + o(1)} \ll p \ll n^{-1/3 + o(1)}$ and we consider this to be an important open problem. It would also be very interesting to prove similarly sharp container theorems for other bipartite graphs $H$.

\section{Containers for multicoloured structures}
\label{sec:many-colours}

All of the problems that we have discussed so far, and many others, are naturally expressed as questions about independent sets in various hypergraphs. There are, however, questions of a very similar flavour that are not easily described in this way but are still amenable to the container method. As an example, consider the problem of enumerating large graphs with no \emph{induced} copy of a given graph $H$. We shall say that a graph $G$ is \emph{induced-$H$-free} if no subset of vertices of $G$ induces a subgraph isomorphic to $H$. As it turns out, it is beneficial to think of an $n$-vertex graph $G$ as the characteristic function of its edge set. A function $g \colon E(K_n) \to \{0, 1\}$ is the characteristic function of an induced-$H$-free graph if and only if for every set $W$ of $v(H)$ vertices of $K_n$, the restriction of $g$ to the set of pairs of vertices of $W$ is not the characteristic function of the edge set of $H$. In particular, viewing $g$ as the set of pairs $\big\{(e, g(e)) \scolon e \in E(K_n)\big\}$, we see that if $g$ represents an induced-$H$-free graph, then it is an independent set in the $\binom{v(H)}{2}$-uniform hypergraph $\cH$ with vertex set $E(K_n) \times \{0, 1\}$ whose edges are the characteristic functions of all copies of $H$ in $K_n$; formally, for every injection $\varphi \colon V(H) \to V(K_n)$, the set
\[
  \Big\{\big(\varphi(u)\varphi(v), 1\big) \scolon uv \in E(H)\Big\} \cup \Big\{\big(\varphi(u)\varphi(v), 0\big) \scolon uv \notin E(H)\Big\}
\]
is an edge of $\cH$. Even though the converse statement is not true and not every independent set of $\cH$ corresponds to an induced-$H$-free graph, since we are usually interested in bounding the number of such graphs from above, the above representation can be useful. In particular, Saxton and Thomason~\cite{ST} applied the container method to the hypergraph $\cH$ described above to reprove the following analogue of~\eqref{eq:EFR}, which was originally obtained by Alekseev~\cite{Al92} and by Bollob\'as and Thomason~\cite{BoTh95, BoTh97}. Letting $\Fnind(H)$ denote the family of all induced-$H$-free graphs with vertex set $\{1, \dotsc, n\}$, we have
\[
  |\Fnind(H)| = 2^{(1-1/\col(H)) \binom{n}{2} + o(n^2)},
\]
where $\col(H)$ is the so-called \emph{colouring number}\footnote{The \emph{colouring number} of a graph $H$ is the largest integer $r$ such that for some pair $(r_1, r_2)$ satisfying $r_1 + r_2 = r$, the vertex set of $H$ cannot be partitioned into $r_1$ cliques and $r_2$ independent sets.} of $H$.

This idea of embedding non-monotone properties (such as the family of induced-$H$-free graphs) into the family of independent sets of an auxiliary hypergraph has been used in several other works. In particular, K\"uhn, Osthus, Townsend, and Zhao~\cite{KuOsToZh17} used it to describe the typical structure of oriented graphs without a transitive tournament of a given order. The recent independent works of Falgas-Ravry, O'Connell, Str\"omberg, and Uzzell~\cite{FROCStUz} and of Terry~\cite{Te} have developed a general framework for studying various enumeration problems in the setting of multicoloured graphs~\cite{FROCStUz} and, more generally, in the very abstract setting of finite (model theoretic) structures~\cite{Te}. In order to illustrate some of the ideas involved in applications of this kind, we will discuss the problem of counting finite metric spaces with bounded integral distances.

\subsection{Counting metric spaces}

Let $\M_n^M$ denote the family of metric spaces on a given set of $n$ points with distances in the set $\{1,\ldots,M\}$. Thus $\M_n^M$ may be viewed as the set of all functions $d \colon E(K_n) \to \{1, \dotsc, M\}$ that satisfy the triangle inequality $d(uv) \le d(uw) + d(wv)$ for all $u, v, w$. Since $x \le y + z$ for all $x, y, z \in \{\lceil M/2 \rceil, \dotsc, M\}$, we have
\begin{equation}
  \label{eq:MnM-lower}
  \big|\M_n^M\big| \ge \left|\left\{\left\lceil \frac{M}{2} \right\rceil, \dotsc, M\right\}\right|^{\binom{n}{2}} = \left\lceil \frac{M+1}{2} \right\rceil^{\binom{n}{2}}.
\end{equation}
Inspired by a continuous version of the model suggested Benjamini (and first studied in~\cite{KoMePeSa}), Mubayi and Terry~\cite{MuTe} proved that for every fixed even $M$, the converse of~\eqref{eq:MnM-lower} holds asymptotically, that is, $|\M_n^M| \le \big( 1 + o(1) \big) \big\lceil \frac{M+1}{2} \big\rceil^{\binom{n}{2}}$ as $n \to \infty$. The problem becomes much more difficult, however, when one allows $M$ to grow with $n$. For example, if $M \gg \sqrt{n}$ then the lower bound
\[
  \big|\M_n^M\big| \ge \left[\left(\frac{1}{2} + \frac{c}{\sqrt{n}} \right) M \right]^{\binom{n}{2}}
\]
for some absolute constant $c > 0$, proved in~\cite{KoMePeSa}, is stronger than~\eqref{eq:MnM-lower}. Balogh and Wagner~\cite{BaWa16} proved strong upper bounds on $|\M_n^M|$ under the assumption that $M \ll n^{1/3} / (\log n)^{4/3+o(1)}$. The following almost optimal estimate was subsequently obtained by Kozma, Meyerovitch, Peled, and Samotij~\cite{KoMePeSa}.

\begin{thm}\label{thm:KMPS}
  There exists a constant $C$ such that
  \begin{equation}\label{eq:thm:KMPS}
    \big|\M_n^M\big| \le \left[\left(\frac{1}{2} + \frac{2}{M} + \frac{C}{\sqrt{n}} \right) M \right]^{\binom{n}{2}}
  \end{equation}
  for all $n$ and $M$.
\end{thm}

Here, we present an argument due to Morris and Samotij that derives a mildly weaker estimate using the hypergraph container lemma. Let $\cH$ be the 3-uniform hypergraph with vertex set $E(K_n) \times \{1, \dotsc, M\}$ whose edges are all triples $\big\{(e_1, d_1), (e_2, d_2), (e_3, d_3)\big\}$ such that $e_1, e_2, e_3$ form a triangle in $K_n$ but $d_{\sigma(1)} + d_{\sigma(2)} < d_{\sigma(3)}$ for some permutation $\sigma$ of $\{1, 2, 3\}$. The crucial observation, already made in~\cite{BaWa16}, is that every metric space $d \colon E(K_n) \to \{1, \dotsc, M\}$, viewed as the set of pairs $\big\{(e, d(e)) \scolon e \in E(K_n)\big\}$, is an independent set of $\cH$. This enables the use of the hypergraph container method for bounding $\big|\M_n^M\big|$ from above. Define the volume of a set $A \subset E(K_n) \times \{1, \dotsc, M\}$, denoted by $\vol(A)$, by
\[
  \vol(A) = \prod_{e \in E(K_n)} \Big|\Big\{ d \in \big\{1, \dotsc, M \big\} \scolon (e, d) \in A \Big\}\Big|
\]
and observe that $A$ contains at most $\vol(A)$ elements of $\M_n^M$. The following supersaturation lemma was proved by Morris and Samotij.

\begin{lemma}\label{lemma:metric-space-supsat}
  Let $n \ge 3$ and $M \ge 1$ be integers and suppose that $A \subset E(K_n) \times \{1, \dotsc, M\}$ satisfies 
  \[
    \vol(A) = \bigg[ \bigg( \frac{1}{2} + \eps \bigg) M \bigg]^{\binom{n}{2}}
  \]
  for some $\eps \ge 10/M$. Then there exist $m \le M$ and a set $A' \subset A$ with $|A'| \le mn^2$, such that the hypergraph $\cH' = \cH[A']$ satisfies 
  \[
    e(\cH') \ge \frac{\eps m^2 M}{50 \log M} \binom{n}{3}, \qquad \Delta_1(\cH') \le 4m^2n, \qquad \text{and} \qquad \Delta_2(\cH') \le 2m.
  \]
\end{lemma}

It is not hard to verify that the hypergraph $\cH'$ given by Lemma~\ref{lemma:metric-space-supsat} satisfies the assumptions of the hypergraph container lemma stated in Section~\ref{keylemma:sec} with $r = \eps n^2 M / (2^{11} \log M)$ and $b = O(n^{3/2})$. Consequently, there exist an absolute constant $\delta$ and a collection $\C$ of subsets of $A'$, with
$$|\C| \le \exp\Big( O\big( n^{3/2}\log(nM) \big) \Big),$$
such that, setting $A_C = C \cup (A \setminus A') = A \setminus (A' \setminus C)$ for each $C \in \C$, the following properties hold:
\begin{enumerate}
\item[$(a)$] every metric space in $A$, viewed as a subset of $E(K_n) \times \{1, \dotsc, M\}$, is contained in $A_C$ for some $C \in \C$, and
\item[$(b)$] $|C| \le |A'| - \delta r$ for every $C \in \C$.
\end{enumerate}
Observe that
\begin{align*}
\vol(A_C) & \, \le \, \left(\frac{M-1}{M}\right)^{|A' \setminus C|} \vol(A) \, \le \, e^{-\delta r / M} \, \vol(A) \\
& \, \le \, e^{-\delta\eps n^2/(2^{11} \log M)} \, \vol(A) \, \le \, \left[\left(\frac{1}{2}+\left(1 - \frac{\delta}{2^{12}\log M}\right)\eps\right)M\right]^{\binom{n}{2}}.
\end{align*}
Since every metric space in $\M_n^M$ is contained in $E(K_n) \times \{1, \dotsc, M\}$, by recursively applying this `breaking down' process to depth $O(\log M)^2$, we obtain a family of 
$$\exp\Big( O\big( n^{3/2} (\log M)^2 \log (nM) \big) \Big)$$ 
subsets of $E(K_n) \times \{1, \dotsc, M\}$, each of volume at most $\big(M/2+10\big)^{\binom{n}{2}}$, that cover all of $\M_n^M$. This implies that
\[
  \big|\M_n^M\big| \le \left[\left(\frac{1}{2} + \frac{10}{M} + \frac{C(\log M)^2\log(nM)}{\sqrt{n}} \right) M \right]^{\binom{n}{2}},
\]
which, as promised, is only slightly weaker than~\eqref{eq:thm:KMPS}. 


\section{An asymmetric container lemma}\label{sec:asymm-cont-lemma}

The approach to studying the family of induced-$H$-free graphs described in the previous section has one (rather subtle) drawback: it embeds $\Fnind(H)$ into the family of independent sets of a $\binom{v(H)}{2}$-uniform hypergraph with $\Theta(n^2)$ vertices. As a result, the hypergraph container lemma produces fingerprints of the same size as for the family of graphs without a clique on $v(H)$ vertices. This precludes the study of various threshold phenomena in the context of sparse induced-$H$-free graphs with the use of the hypergraph container lemma presented in Section~\ref{keylemma:sec}; this is in sharp contrast with the non-induced case, where the container method proved very useful.

In order to alleviate this shortcoming, Morris, Samotij, and Saxton~\cite{MSS} proved a version of the hypergraph container lemma for $2$-coloured structures that takes into account the possible asymmetries between the two colours. We shall not give the precise statement of this new container lemma here (since it is rather technical), but we would like to emphasize the following key fact: it enables one to construct families of containers for induced-$H$-free graphs with fingerprints of size $\Theta(n^{2-1/m_2(H)})$, as in the non-induced case. 

To demonstrate the power of the asymmetric container lemma, the following application was given in~\cite{MSS}. Let us say that a graph $G$ is \emph{$\eps$-close to a split graph} if there exists a partition $V(G) = A \cup B$ such that $e(G[A]) \ge (1-\eps)\binom{|A|}{2}$ and $e(G[B]) \le \eps e(G)$.

\begin{thm}\label{thm:MSS}
  For every $\eps > 0$, there exists a $\delta > 0$ such that the following holds. Let $G$ be a~uniformly chosen induced-$C_4$-free graph with vertex set $\{1, \dotsc, n\}$ and $m$ edges.
  \begin{enumerate}
  \item[$(a)$] If $n \ll m \ll \delta n^{4/3} (\log n)^{1/3}$, then a.a.s.\ $G$ is not $1/4$-close to a split graph.
  \item[$(b)$] If $n^{4/3} (\log n)^4 \le m \le \delta n^2$, then a.a.s.\ $G$ is $\eps$-close to a split graph.
  \end{enumerate}
\end{thm}

Theorem~\ref{thm:MSS} has the following interesting consequence: it allows one to determine the number of edges in (and, sometimes, also the typical structure of) the binomial random graph $G(n,p)$ conditioned on the event that it does not contain an induced copy of $C_4$. Let us denote by $\Gnpind(C_4)$ the random graph chosen according to this conditional distribution. 

\begin{cor}\label{cor:GnpC4free}
The following bounds hold asymptotically almost surely as $n \to \infty$:
  \[
    e \big( \Gnpind(C_4) \big) \, = \, 
    \begin{cases}
      \, \big( 1 + o(1) \big) p {n \choose 2} & \textup{if $n^{-1} \ll p \ll n^{-2/3}$}, \\
      \, n^{4/3} (\log n)^{O(1)} & \textup{if $n^{-2/3} \le p \le n^{-1/3} (\log n)^4$}, \\
      \, \Theta\big( p^2 n^2 / \log(1/p) \big) & \textup{if $p \ge n^{-1/3} (\log n)^4$}.
    \end{cases}
  \]
\end{cor}

We would like to emphasize the (surprising) similarity between the statements of Theorem~\ref{thm:randomturan} and Corollary~\ref{cor:GnpC4free}. In particular, the graph of $p \mapsto e \big( \Gnpind(C_4) \big)$ contains exactly the same `long flat segment' as the graph of $p \mapsto \ex \big( G(n,p), C_4 \big)$, even though the shape of the two graphs above this range is quite different. We do not yet fully understand this phenomenon and it would be interesting to investigate whether or not the function $p \mapsto e \big( \Gnpind(H) \big)$ exhibits similar behaviour for other bipartite graphs $H$.

\section{Hypergraphs of unbounded uniformity}\label{unbounded:sec}


Since the hypergraph container lemma provides explicit dependencies between the various parameters in its statement, it is possible to apply the container method even when the uniformity of the hypergraph considered is a growing function of the number of its vertices. Perhaps the first result of this flavour was obtained by Mousset, Nenadov, and Steger~\cite{MoNeSt14}, who proved an upper bound of $2^{\ex(n, K_r) + o(n^2/r)}$ on the number of $n$-vertex $K_r$-free graphs for all $r \le (\log_2n)^{1/4}/2$. Subsequently, Balogh, Bushaw, Collares, Liu, Morris, and Sharifzadeh~\cite{BBCLMS} strengthened this result by establishing the following precise description of the typical structure of large $K_r$-free graphs.

\begin{thm}
  If $r \le (\log_2 n)^{1/4}$, then almost all $K_r$-free graphs with $n$ vertices are $(r-1)$-partite.
\end{thm}

Around the same time, the container method applied to hypergraphs with unbounded uniformity was used to analyse Ramsey properties of random graphs and hypergraphs, leading to improved upper bounds on several well-studied functions. In particular, R\"odl, Ruci\'nski, and Schacht~\cite{RoRuSc17} gave the following upper bound on the so-called Folkman numbers. 

\begin{thm}
For all integers $k \ge 3$ and $r \ge 2$, there exists a $K_{k+1}$-free graph with
$$\exp\big( Ck^4\log k + k^3r\log r \big)$$ 
vertices, such that every $r$-colouring of its edges contains a monochromatic copy of $K_k$. 
\end{thm}

The previously best known bound was doubly exponential in $k$, even in the case $r = 2$. Not long afterwards, Conlon, Dellamonica, La Fleur, R\"odl, and Schacht~\cite{CoDeLFRoSc17} used a similar method to prove the following strong upper bounds on the induced Ramsey numbers of hypergraphs. Define the tower functions $t_k(x)$ by $t_1(x) = x$ and $t_{i+1}(x) = 2^{t_i(x)}$ for each $i \ge 1$. 

\begin{thm}
For each $k \ge 3$ and $r \ge 2$, there exists $c$ such that the following holds. For every $k$-uniform hypergraph $F$ on $m$ vertices, there exists a $k$-uniform hypergraph $G$ on~$t_k(cm)$ vertices, such that every $r$-colouring of $E(G)$ contains a monochromatic induced copy of~$F$.  
\end{thm}


Finally, let us mention a recent result of Balogh and Solymosi~\cite{BS}, whose proof is similar to that of Theorem~\ref{thm:geometry}, which we outlined in Section~\ref{geometry:sec}. Given a family $\cF$ of subsets of an $n$-element set $\Omega$, an \emph{$\eps$-net} of $\cF$ is a set $A \subset \Omega$ that intersects every member of $\cF$ with at least $\eps n$ elements. The concept of an $\eps$-net plays an important role in computer science, for example in computational geometry and approximation theory. In a seminal paper, Haussler and Welzl~\cite{HaWe87} proved that every set system with VC-dimension\footnote{The VC-dimension (VC stands for Vapnik--Chervonenkis) of a family $\cF$ of subsets of $\Omega$ is the largest size of a~set $X \subset \Omega$ such that the set $\{A \cap X \scolon A \in \cF\}$ has $2^{|X|}$ elements.} $d$ has an $\eps$-net of size $O\big( (d/\eps) \log(d/\eps) \big)$. It was believed for more than twenty years that for `geometric' families, the $\log(d/\eps)$ factor can be removed; however, this was disproved by Alon~\cite{A12}, who constructed, for each $C > 0$, a set of points in the plane such that the smallest $\eps$-net for the family of lines (whose VC-dimension is $2$) has size at least~$C/\eps$. 

By applying the container method to the hypergraph of collinear $k$-tuples in the $k$-dimensional $2^{k^2} \times \cdots \times 2^{k^2}$ integer grid, Balogh and Solymosi~\cite{BS} gave the following stronger lower bound. 

\begin{thm}\label{thm:epsnets}
  For each $\eps > 0$, there exists a set $S \subset \R^2$ such that the following holds. If $T \subset S$ intersects every line that contains at least $\eps |S|$ elements of $S$, then
  \[
    |T| \ge \frac{1}{\eps} \left( \log \frac{1}{\eps} \right)^{1/3 + o(1)}.
  \]
\end{thm}

It was conjectured by Alon~\cite{A12} that there are sets of points in the plane whose smallest $\eps$-nets (for the family of lines) contain $\Omega\big(1/\eps \log(1/\eps)\big)$ points.

\section{Some further applications}\label{more:applications:sec}

There are numerous applications of the method of containers that we do not have space to discuss in detail. Still, we would like to finish this survey by briefly mentioning just a few of them.

\subsection{List colouring}

A hypergraph $\cH$ is said to be \emph{$k$-choosable} if for every assignment of a list $L_v$ of $k$ colours to each vertex $v$ of $\cH$, it is possible to choose for each $v$ a colour from the list $L_v$ in such a way that no edge of $\cH$ has all its vertices of the same colour. The smallest $k$ for which $\cH$ is $k$-choosable is usually called the \emph{list chromatic number} of $\cH$ and denoted by $\chi_\ell(\cH)$. Alon~\cite{Al93,Al00} showed that for graphs, the list chromatic number grows with the minimum degree, in stark contrast with the usual chromatic number; more precisely, $\chi_\ell(G) \ge \big( 1/2 + o(1) \big) \log_2 \delta(G)$ for every graph $G$. The following generalisation of this result, which also improves the constant $1/2$, was proved by Saxton and Thomason~\cite{ST}, see also~\cite{SaTh12,ST16}.

\begin{thm}
  Let $\cH$ be a $k$-uniform hypergraph with average degree $d$ and $\Delta_2(\cH) = 1$. Then, as $d \to \infty$,
  \[
    \chi_\ell(\cH) \ge \left(\frac{1}{(k-1)^2} + o(1)\right) \log_k d.
  \]
Moreover, if $\cH$ is $d$-regular, then
  \[
    \chi_\ell(\cH) \ge \left(\frac{1}{k-1} + o(1)\right) \log_kd.
  \]
\end{thm}

We remark that proving lower bounds for the list chromatic number of simple hypergraphs was one of the original motivations driving the development of the method of hypergraph containers.

\subsection{Additive combinatorics}

The method of hypergraph containers has been applied to a number of different number-theoretic objects, including sum-free sets~\cite{ABMS1,ABMS2,BLST2,BLST1}, Sidon sets~\cite{ST16}, sets containing no $k$-term arithmetic progression~\cite{BLS,BMS}, and general systems of linear equations~\cite{ST16}. (See also~\cite{Gr04,GrRu04,Sap03} for early applications of the container method to sum-free sets and~\cite{DeKoLeRoSa-B3, DeKoLeRoSa-Bh-16, DeKoLeRoSa-Bh, KoLeRoSa} for applications of graph containers to $B_h$-sets.) Here we will mention just three of these results. 

Let us begin by recalling that a {\it Sidon set} is a set of integers containing no non-trivial solutions of the equation $x+y=z+w$. Results of Chowla, Erd\H os, Singer, and Tur\'an from the 1940s imply that the maximum size of a Sidon set in $\{1,\ldots,n\}$ is $\big( 1 + o(1) \big) \sqrt n$ and it was conjectured by Cameron and Erd{\H{o}}s~\cite{CamErdos} that the number of such sets is $2^{(1+o(1))\sqrt n}$. This conjecture was disproved by Saxton and Thomason~\cite{ST16}, who gave a construction of $2^{(1+\eps) \sqrt n}$ Sidon sets (for some $\eps > 0$), and also used the hypergraph container method to reprove the following theorem, which was originally obtained in~\cite{KoLeRoSa} using the graph container method.

\begin{thm}\label{sidon}
There are $2^{O(\sqrt n)}$ Sidon sets in $\{1,\ldots,n\}$. 
\end{thm}

Dellamonica, Kohayakawa, Lee, R\"odl, and Samotij~\cite{DeKoLeRoSa-Bh} later generalized these results to $B_h$-sets, that is, set of integers containing no non-trivial solutions of the equation $x_1+ \ldots + x_h = y_1 + \ldots + y_h$.

The second result we would like to state was proved by Balogh, Liu, and Sharifzadeh~\cite{BLS}, and inspired the proof presented in Section~\ref{counting:Hfree:sec}. Let $r_k(n)$ be the largest size of a subset of $\{1, \dotsc, n\}$ containing no $k$-term arithmetic progressions.

\begin{thm}\label{thm:BLS}
For each integer $k \ge 3$, there exist a constant $C$ and infinitely many $n \in \N$ such that there are at most $2^{C r_k(n)}$ subsets of $\{1, \dotsc, n\}$ containing no $k$-term arithmetic progression.
\end{thm}

We recall (see, e.g.,~\cite{GowersErdos}) that obtaining good bounds on $r_k(n)$ is a well-studied and notoriously difficult problem. The proof of Theorem~\ref{thm:BLS} avoids these difficulties by exploiting merely the `self-similarity' property of the hypergraph encoding arithmetic progressions in $\{1,\ldots,n\}$, cf.~the discussion in Section~\ref{keylemma:sec} and the proof of Lemma~\ref{lem:C4:super}.

The final result we would like to mention was one of the first applications of (and original motivations for the development of) the method of hypergraph containers. Recall that the Cameron--Erd\H{o}s conjecture, proved by Green~\cite{Gr04} and, independently, by Sapozhenko~\cite{Sap03}, states that there are only $O(2^{n/2})$ sum-free subsets of $\{1,\ldots,n\}$. The following sparse analogue of the Cameron--Erd\H{o}s  conjecture was proved by Alon, Balogh, Morris, and Samotij~\cite{ABMS2} using an early version of the hypergraph container lemma for 3-uniform hypergraphs.

\begin{thm}\label{CEthm}
There exists a constant $C$ such that, for every $n \in \N$ and every $1 \le m \le \lceil n/2 \rceil$, the set $\{1,\ldots,n\}$ contains at most $2^{Cn/m} {\lceil n/2 \rceil \choose m}$ sum-free sets of size $m$.
\end{thm}

We remark that if $m \ge \sqrt{n}$, then Theorem~\ref{CEthm} is sharp up to the value of $C$, since in this case there is a constant $c > 0$ such that there are at least $2^{c n/m} {n/2 \choose m}$ sum-free $m$-subsets of $\{1,\ldots,n\}$. For smaller values of $m$ the answer is different, but the problem in that range is much easier and can be solved using standard techniques. Let us also mention that in the case $m \gg \sqrt{n \log n}$, the structure of a typical sum-free $m$-subset of $\{1,\ldots,n\}$ was also determined quite precisely in~\cite{ABMS2}.

Finally, we would like to note that, although the statements of Theorems~\ref{sidon},~\ref{thm:BLS} and~\ref{CEthm} are somewhat similar, the difficulties encountered during their proofs are completely different.

\subsection{Sharp thresholds for Ramsey properties}

Given an integer $k \ge 3$, let us say that a set $A \subset \Z_n$ has the \emph{van der Waerden property} for $k$ if every 2-colouring of the elements of $A$ contains a monochromatic $k$-term arithmetic progression; denote this by $A \to (k\text{-AP})$. R\"odl and Ruci\'nski~\cite{RoRu95} determined the threshold for the van der Waerden property in random subsets of $\Z_n$ for every~$k \in \N$. Combining the sharp threshold technology of Friedgut~\cite{F99} with the method of hypergraph containers, Friedgut, H\'an, Person, and Schacht~\cite{FrHaPeSc16} proved that this threshold is sharp. Let us write $\Z_{n,p}$ to denote a $p$-random subset of $\Z_n$ (i.e., each element is included independently with probability $p$). 

\begin{thm}
For every $k \ge 3$, there exist constants $c_1 > c_0 > 0$ and a function $p_c \colon \N \to [0,1]$ satisfying $c_0 n^{-1/(k-1)} < p_c(n) < c_1 n^{-1/(k-1)}$ for every $n \in \N$, such that, for every $\eps > 0$,
  \[
    \Pr\Big( \Z_{n,p} \to  \big( \text{$k$-AP} \big) \Big) \, \to \,
    \begin{cases}
     \, 0 & \text{if } \, p \le (1-\eps) \, p_c(n),\\[+0.05ex]
    \,  1 & \text{if } \, p \ge (1+\eps) \, p_c(n),
    \end{cases}
  \]
  as $n \to \infty$. 
\end{thm}

The existence of a sharp threshold in the context of Ramsey's theorem for the triangle was obtained several years earlier, by Friedgut, R\"odl, Ruci\'nski, and Tetali~\cite{FRRT}. Very recently, using similar methods to those in~\cite{FrHaPeSc16}, Schacht and Schulenburg~\cite{SchSch} gave a simpler proof of this theorem and also generalised it to a large family of graphs, including all odd cycles.

\subsection{Maximal triangle-free graphs and sum-free sets}

In contrast to the large body of work devoted to counting and describing the typical structure of $H$-free graphs, relatively little is known about $H$-free graphs that are \emph{maximal} (with respect to the subgraph relation). The following construction shows that there are at least $2^{n^2/8}$ maximal triangle-free graphs with vertex set $\{1, \dotsc, n\}$. Fix a partition $X \cup Y = \{1, \dotsc, n\}$ with $|X|$ even. Define $G$ by letting $G[X]$ be a perfect matching, leaving $G[Y]$ empty, and adding to $E(G)$ exactly one of $xy$ or $x'y$ for every edge $xx' \in E(G[X])$ and every $y \in Y$. It is easy to verify that all such graphs are triangle-free and that almost all of them are maximal. 

Using the container theorem for triangle-free graphs (Theorem~\ref{thm:CT:triangles}), Balogh and Pet\v{r}\'\i\v{c}kov\'a~\cite{BaPe14} proved that the construction above is close to optimal by showing that there are at most $2^{n^2/8 + o(n^2)}$ maximal triangle-free graphs on $\{1, \dotsc, n\}$. Following this breakthrough, Balogh, Liu, Pet\v{r}\'\i\v{c}kov\'a, and Sharifzadeh~\cite{BaLiPeSh15} proved the following much stronger theorem, which states that in fact almost all maximal triangle-free graphs can be constructed in this way.

\begin{thm}
  For almost every maximal triangle-free graph $G$ on $\{1, \dotsc, n\}$, there is a vertex partition $X \cup Y$ such that $G[X]$ is a perfect matching and $Y$ is an independent set.
\end{thm}

A similar result for sum-free sets was obtained by Balogh, Liu, Sharifzadeh, and Treglown~\cite{BLST2,BLST1}, who determined the number of maximal sum-free subsets of $\{1, \dotsc, n\}$ asymptotically. However, the problem of estimating the number of maximal $H$-free graphs for a general graph $H$ is still wide open. In particular, generalizing the results of~\cite{BaLiPeSh15,BaPe14} to the family of maximal $K_k$-free graphs seems to be a very interesting and difficult open problem. 

\subsection{Containers for rooted hypergraphs} 

A family $\cF$ of finite sets is \emph{union-free} if $A \cup B \ne C$ for every three distinct sets $A, B, C \in \cF$. Kleitman~\cite{K76} proved that every union-free family in $\{1,\dotsc,n\}$ contains at most $\big( 1 + o(1) \big) \binom{n}{n/2}$ sets; this is best possible as the family of all $\lfloor n/2 \rfloor$-element subsets of $\{1, \dotsc, n\}$ is union-free. Balogh and Wagner~\cite{BaWaUnion} proved the following natural counting counterpart of Kleitman's theorem, confirming a conjecture of Burosch, Demetrovics, Katona, Kleitman, and Sapozhenko~\cite{BDKKS}.

\begin{thm}
There are $2^{( 1 + o(1)) \binom{n}{n/2}}$ union-free families in $\{1,\ldots,n\}$. 
\end{thm}

It is natural to attempt to prove this theorem by applying the container method to the 3-uniform hypergraph $\cH$ that encodes triples $\{A,B,C\}$ with $A \cup B = C$. However, there is a problem: for any pair $(B,C)$, there exist $2^{|B|}$ sets $A$ such that $A \cup B = C$ and this means that $\Delta_2(\cH)$ is too large for a naive application of the hypergraph container lemma. In order to overcome this difficulty, Balogh and Wagner developed in~\cite{BaWaUnion} a new container  theorem for `rooted' hypergraphs (each edge has a designated root vertex) that exploits the asymmetry of the identity $A \cup B = C$. In particular, note that while the degree of a pair $(B, C)$ can be large, the pair $\{A,B\}$ uniquely determines $C$; it turns out that this is sufficient to prove a suitable container theorem. We refer the reader to~\cite{BaWaUnion} for the details.

\subsection{Probabilistic embedding in sparse graphs}

The celebrated regularity lemma of Szemer\'edi~\cite{Sz78} states that, roughly speaking, the vertex set of every graph can be divided into a bounded number of parts in such a way that most of the bipartite subgraphs induced by pairs of parts are pseudorandom; such a partition is called a~\emph{regular partition}. The strength of the regularity lemma stems from the so-called counting and embedding lemmas, which tell us approximately how many copies of a particular subgraph a graph $G$ contains in terms of basic parameters of the regular partition of $G$. While the original statement of the regularity lemma applied only to dense graphs (i.e., $n$-vertex graphs with $\Omega(n^2)$ edges), the works of Kohayakawa~\cite{Ko97}, R\"odl~(unpublished), and Scott~\cite{Sc11} provide extensions of the lemma that are applicable to sparse graphs. However, these extensions come with a major caveat: the counting and embedding lemmas do not extend to sparse graphs; this unfortunate fact was observed by \L uczak. Nevertheless, it seemed likely that such atypical graphs that fail the counting or embedding lemmas are so rare that they typically do not appear in random graphs. This belief was formalised in a conjecture of Kohayakawa, \L uczak, and R\"odl~\cite{KLR2}, which can be seen as a `probabilistic' version of the embedding lemma. 

The proof of this conjecture, discovered by the authors of this survey~\cite{BMS} and by Saxton and Thomason~\cite{ST}, was one of the original applications of the hypergraph container lemma. Let us mention here that a closely related result was proved around the same time by Conlon, Gowers, Samotij, and Schacht~\cite{CoGoSaSc14}. A~strengthening of the K\L R conjecture, a `probabilistic' version of the counting lemma, proposed by Gerke, Marciniszyn, and Steger~\cite{GeMaSt07}, remains open.

\section*{Acknowledgements}

The authors would like to thank Noga Alon, B\'ela Bollob\'as, David Conlon, Yoshi Kohayakawa, Gady Kozma, Tom Meyerovitch, Ron Peled, David Saxton, and Andrew Thomason for many stimulating discussions over the course of the last several years. These discussions have had a significant impact on the development of the container method. Last but not least, we would like to again acknowledge David Conlon, Ehud Friedgut, Tim Gowers, Daniel Kleitman, Vojta R\"odl, Mathias Schacht, and Kenneth Wilson, whose work on enumeration of $C_4$-free graphs and on extremal and Ramsey properties of random discrete structures inspired and influenced our investigations that led to the hypergraph container theorems.

\bibliographystyle{amsplain}
\bibliography{ICM-containers}

\end{document}